\documentclass[12pt]{article}
\usepackage{epic,eepic}
\usepackage[dvips]{epsfig}
\usepackage{color}
\usepackage{amsmath, amscd, amssymb}
\usepackage{amssymb}
\def\draw #1 by #2 (#3){
  \vbox to #2{
    \hrule width #1 height 0pt depth 0pt
    \vfill
    \special{picture #3}
    }
  }
\def\scaleddraw #1 by #2 (#3 scaled #4){{
  \dimen0=#1 \dimen1=#2
  \divide\dimen0 by 1000 \multiply\dimen0 by #4
  \divide\dimen1 by 1000 \multiply\dimen1 by #4
  \draw \dimen0 by \dimen1 (#3 scaled #4)}
  }

\begin{document}
\renewcommand{\labelenumi}{\theenumi}
\newcommand{\qed}{\mbox{\raisebox{0.7ex}{\fbox{}}}}
\newtheorem{theorem}{Theorem}
\newtheorem{example}{Example}
\newtheorem{problem}[theorem]{Problem}
\newtheorem{defin}[theorem]{Definition}
\newtheorem{lemma}[theorem]{Lemma}
\newtheorem{corollary}[theorem]{Corollary}
\newtheorem{nt}{Note}
\newtheorem{proposition}[theorem]{Proposition}
\renewcommand{\thent}{}
\newenvironment{pf}{\medskip\noindent{\textbf{Proof}:  \hspace*{-.4cm}}\enspace}{\hfill \qed \medskip \newline}
\newenvironment{defn}{\begin{defin}}{\end{defin}}{\vspace{-0.5cm}}
\newenvironment{lem}{\begin{lemma}}{\end{lemma}}{\vspace{-0.5cm}}
\newenvironment{cor}{\begin{corollary}}{\end{corollary}}{\vspace{-0.5cm}}
\newenvironment{thm}{\begin{theorem} }{\end{theorem}}{\vspace{-0.5cm}}
\newenvironment{pbm}{\begin{problem} \em}{\end{problem}}{\vspace{-0.5cm}}
\newenvironment{note}{\begin{nt} \em}{\end{nt}}{\vspace{-0.5cm}}
\newenvironment{exa}{\begin{example} \em}{\end{example}}{\vspace{-0.5cm}}
\newenvironment{pro}{\begin{proposition} \em}{\end{proposition}}{\vspace{-0.5cm}}
\setlength{\unitlength}{12pt}
\newcommand{\comb}[2]{\mbox{$\left(\!\!\begin{array}{c}
            {#1} \\[-0.5ex] {#2} \end{array}\!\!\right)$}}
\renewcommand{\labelenumi}{(\theenumi)}
\renewcommand{\b}{\beta}
\newcounter{myfig}
\newcounter{mytab}
\def\mod{\hbox{\rm mod }}
\def\scaleddraw #1 by #2 (#3 scaled #4){{
  \dimen0=#1 \dimen1=#2
  \divide\dimen0 by 1000 \multiply\dimen0 by #4
  \divide\dimen1 by 1000 \multiply\dimen1 by #4
  \draw \dimen0 by \dimen1 (#3 scaled #4)}
  }
\newcommand{\Aut}{\mbox{\rm Aut}}
\newcommand{\w}{\omega}
\def\r{\rho}
\newcommand{\DbF}{D \times^{\phi} F}
\newcommand{\autF}{{\tiny\Aut{\scriptscriptstyle(\!F\!)}}}
\def\Cay{\mbox{\rm Cay}}
\def\a{\alpha}
\newcommand{\C}[1]{\mathcal #1}
\newcommand{\B}[1]{\mathbb #1}
\newcommand{\F}[1]{\mathfrak #1}
\title{\textbf{Shilla distance-regular graphs}}
  \author{Jack H. Koolen, Jongyook Park \\
\
{\small Department of Mathematics, POSTECH, Pohang, 790-784, Korea }\\
          {\small E-mail: koolen@postech.ac.kr, jongyook@postech.ac.kr}\\
}
\date{}
\maketitle

\begin{abstract}

A Shilla distance-regular graph $\Gamma$ (say with valency $k$) is a distance-regular graph with diameter 3 such that its second largest eigenvalue equals to $a_3$. We will show that $a_3$ divides $k$ for a Shilla distance-regular graph $\Gamma$, and for $\Gamma$ we define $b=b(\Gamma):=\frac{k}{a_3}$. In this paper we will show that there are finitely many Shilla distance-regular graphs $\Gamma$ with fixed $b(\Gamma)\geq 2$. Also, we will classify  Shilla distance-regular graphs with $b(\Gamma)=2$ and $b(\Gamma)=3$. Furthermore, we will give a new existence condition for distance-regular graphs, in general.

\bigskip
\noindent
 {{\bf Key Words: distance-regular graph; Existence condition; Terwilliger graph } }
\\
\noindent
 {{\bf 2000 Mathematics Subject Classification: 05E30} }
\\
\end{abstract}

\section{Introduction}

~~~~~~In this paper we study distance-regular graphs $\Gamma$ with diameter 3. (For definitions, see next section.) For a distance-regular graph with diameter 3, we will show that the second largest eigenvalue $\theta_1$  is at least  min$\{\frac{a_{1}+\sqrt{a_{1}^{2}+4k}}{2}, a_{3}\}$, where $k$ is the valency (see, Lemma~\ref{10} below), and that $\theta_1=a_3$ if and only if $\theta_1$=$\frac{a_{1}+\sqrt{a_{1}^{2}+4k}}{2}$. A distance-regular graph $\Gamma$ with diameter 3 is called {\em Shilla} if $\theta_1=a_3$. It follows that for a Shilla distance-regular graph $\Gamma$, $a_3$ divides $k$ and we will put $b(\Gamma):=\frac{k}{a_3}$. In this paper we will show that there exist finitely many (non-isomorphic) Shilla distance-regular graphs with fixed $b(\Gamma)\geq 2$. This result relies on a new existence condition, Theorem~\ref{3}, for distance-regular graphs. Furthermore we will classify Shilla distance-regular graphs $\Gamma$ with $b(\Gamma)\in\{2,3\}$.

This paper is organized as follows:\\
In Section 2, we will give definitions. In Section 3, we give the new existence condition for distance-regular graphs, and in Section 4 we will discuss Shilla distance-regular graphs.

\section{Definitions and preliminaries}

~~~~Suppose that $\Gamma$ is a connected graph with the vertex set $V(\Gamma)$ and the edge set $E(\Gamma)$, where $E(\Gamma)$ consists of the unordered pairs of adjacent two vertices. The distance $d_{\Gamma}(x,y)$ between any two vertices $x$, $y$ of $\Gamma$ is the length of a shortest path between $x$ and $y$ in $\Gamma$.\\

Let $\Gamma$ be a connected graph. For a vertex $x\in V(\Gamma)$, define $\Gamma_i(x)$ to be the set of vertices which are at distance precisely $i$ from $x$ ($0\leq i\leq D$) where $D$ := max$\{d_{\Gamma}(x,y) \mid x, y \in V(\Gamma)\}$ is the diameter of $\Gamma$.  In addition, define $\Gamma_{-1}(x):=\emptyset$ and $\Gamma_{D+1}(x) := \emptyset.$
We will write $\Gamma(x)$ instead of $\Gamma_1(x)$ and we denote $x\sim _{\Gamma}y$ or simply $x\sim y$ if two vertices $x$ and $y$ are adjacent in $\Gamma$. For $x_1,x_2,\cdots,x_l\in V(\Gamma)$, define
 $$\Gamma(x_1,x_2,\cdots,x_l):=\displaystyle\bigcap^l_{i=1}\Gamma(x_i) .$$\\

A connected graph $\Gamma$ with diameter $D$ is called {\em{distance-regular}} if there are integers $b_i, c_i$ $(0\leq i\leq D)$ such that for any two vertices $x, y \in V(\Gamma)$ with $d_{\Gamma}(x, y)=i$, there are precisely $c_i$ neighbors of $y$ in $\Gamma_{i-1}(x)$ and $b_i$ neighbors of $y$ in $\Gamma_{i+1}(x)$. In particular, distance-regular graph $\Gamma$ is regular with valency $k := b_0$ and we define $a_i := k-b_i-c_i$ for notational convenience.  Note that $a_i=\mid\Gamma(y)\cap\Gamma_i(x)\mid$ holds for any two vertices $x, y$ with $d_{\Gamma}(x, y)=i$ $(0\leq i\leq D).$ For a distance-regular graph $\Gamma$ and a vertex $x\in V(\Gamma)$, we denote $k_i:=|\Gamma_i(x)|$. The numbers $a_i$, $b_{i-1}$ and $c_i$ $(1\leq i\leq D)$ are called the {\em{intersection~numbers}} of $\Gamma$, and they satisfy the following three conditions:\\

\vspace{-0.3cm}$(i)$ $k=b_0> b_1\geq \cdots \geq b_{D-1}$;\\

$(ii)$ $1=c_1\leq c_2\leq \cdots \leq c_D$;\\

$(iii)$ $b_i\geq c_j$ if $i+j\leq D$.\\

The array $\{b_0,b_1,\cdots,b_{D-1};c_1,c_2,\cdots,c_D\}$ is called the {\em{intersection~array}} of a distance-regular graph $\Gamma$.
Suppose that $\Gamma$ is a distance-regular graph with valency $k \geq 2$ and diameter $D \geq 2$, and let $A_i$ be the matrix of $\Gamma$ such that the rows and the columns of $A_i$ are indexed by $V(\Gamma)$ and the ($x, y$)-entry of $A_i$ equals 1 whenever $d_{\Gamma}(x,y)=i$ and 0 otherwise. We will denote the adjacency matrix of $\Gamma$ as $A$ instead of
$A_1$. Then $\Gamma$ has exactly ($D+1$) distinct eigenvalues, say $k=\theta_0 > \theta_1 > \cdots > \theta_D$, and let $m_i$ be the multiplicity of $\theta_i$ ($0\leq i \leq D$) , where an eigenvalue of $\Gamma$ is that of $A$.\\

For an eigenvalue $\theta$ of $\Gamma$, the sequence $u_0=u_0(\theta)=1$, $u_1=u_1(\theta)=\frac{\theta}{k}$, $u_i=u_i(\theta)$ ($2\leq i \leq D$) satisfying $$c_iu_{i-1}(\theta)+a_iu_i(\theta)+b_iu_{i+1}(\theta)=\theta u_i(\theta)$$ is called the {\em{standard sequence}} corresponding to the eigenvalue $\theta$.\\

N. Biggs\cite[p.131]{bcn} showed that for an eigenvalue $\theta$ of a distance-regular graph $\Gamma$, its multiplicity $m$ is given by

 \begin{equation}\label{e0}
 m=\frac{\mid V(\Gamma)\mid}{\tiny{\displaystyle \sum_{i=0}^Dk_iu_i(\theta)^2}}.
 \end{equation}

The Bose-Mesner algebra $M$ for a distance-regular graph $\Gamma$ is the matrix algebra generated by the adjacency matrix $A$ of $\Gamma$. A basis of $M$ is $\{A_i~ |~ i=0,\cdots,D\}$, where $A_0=I$. The algebra $M$ has also a basis consisting of primitive idempotents $\{E_0=\frac{1}{n}J, E_1,\cdots,E_D\}$, where $n=|V(\Gamma)|$ and $E_i$ is the orthogonal projection onto the eigenspace of $\theta_i$. Under the componentwise multiplication $\circ$,  $E_i\circ E_j=\frac{1}{n}\displaystyle\sum^D_{k=0} q^k_{ij}E_k$.
The numbers $q^k_{ij}$ ($0\leq i,j,k\leq D$) are called the {\em Krein parameters} of $\Gamma$ and are always non-negative by Delsarte~\cite[Theorem 2.3.2]{bcn}. We say that $\Gamma$ is {\em Q-polynomial} if there is an order of the primitive idempotents $E_0=\frac{1}{n}J, E_1,\cdots,E_D$ such that $q^k_{1j}=0$ if $|j-k|>1$. We say that $\Gamma$ is  $Q$-polynomial with respect to $\theta$ if $E_1$ is the orthogonal projection on the eigenspace of $\theta$.\\

In this paper we say that an intersection array is {\em feasible} if it satisfies the following four conditions:\\

\vspace{-0.3cm}\hspace{-0.6cm}$(i)$ all its intersection numbers $p^i_{jl}$ are integral;\\ \hspace{0.3cm}(where $p^i_{jl}=|\{z\mid d_{\Gamma}(x,z)=j, d_{\Gamma}(y,z)=l\}|$~ for any vertices $x$ and $y$ at distance $i$)\\

\hspace{-0.6cm}$(ii)$ all the multiplicities are positive integers;\\

\hspace{-0.6cm}$(iii)$ for any $0\leq i\leq D$, $k_ia_i$ is even;\\

\hspace{-0.6cm}$(iv)$ all Krein parameters are non-negative.\\

Recall that a {\em clique} of a graph is a set of mutually adjacent vertices and that a {\em co-clique} of a graph is a set of vertices with no edges. For a graph $\Gamma$, the {\em local graph} at a vertex $x\in V(\Gamma)$ is the subgraph induced by $\Gamma(x)$ in $\Gamma$ and we denote it by $\Delta(x)$. Let $\Delta$ be a graph. We say $\Gamma$ is locally $\Delta$ if the local graph $\Delta(x)$ is isomorphic to $\Delta$ for all vertices $x\in V(\Gamma)$. An {\em$order (s,t)$-graph} is a graph such that each $\Delta(x)$ is the disjoint union of $t+1$ copies of $(s+1)$-cliques. A {\em Terwilliger graph} is a connected non-complete graph $\Gamma$ such that, for any two vertices $u, v$ at distance two, the subgraph induced by $\Gamma(u, v)$ in $\Gamma$ is a clique of size $\mu$ (for some fixed $\mu \geq 1$).\\

 Recall the following interlacing result.

\begin{theorem}\label{0}\rm\textbf{(cf. Haemers\cite{heam})}   Let $A$ be a real symmetric $n\times n$ matrix and let $B$ be a principal submatrix of $A$ with order $m \times m$. Then, for $i=1,\cdots , m$, $$\theta_{n-m+i}(A)\leq \theta_i(B)\leq \theta_i(A).$$\\
\end{theorem}

\vspace{-1.5cm}
\section{A new existence condition}

~~~~~In this section, we will give a new existence condition, Theorem~\ref{3}, for distance-regular graphs. To do this we first show Lemma~\ref{1} and Proposition~\ref{2}.

 \begin{lemma}\label{1} \rm{Let} $\Gamma$ be a distance-regular graph with valency $k$ and diameter $D \geq 2$. Let $x$ be a vertex of $\Gamma$ and let $\bar{C}$ be a co-clique of size $s\geq 2$ in the local graph $\Delta(x)$ at $x$. Then\\ $$c_2-1 \geq \frac{s(a_1+1)-k}{{s \choose 2}}.$$
\end{lemma}
\begin{pf}
Let $V(\bar{C})=\{y_1, y_2, \cdots, y_s \}$. Since $d_{\Gamma}(y_i, y_j)$=2,  $\mid\Gamma(x,y_i,y_j)\mid$ $\leq c_2-1$ holds for any $i\neq j$. Then by the principle of inclusion and exclusion,\\

\hspace{3cm}$k=\mid\Gamma(x)\mid\geq| \displaystyle \bigcup_{i=1}^s(\Gamma(x,y_i)\cup\{y_i\})|$\

 \hspace{4.95cm}$ \geq \displaystyle \sum^s_{i=1}\mid \Gamma(x,y_i)\cup\{y_i\}\mid - \displaystyle \sum_{1 \leq i < j \leq s}\mid\Gamma(x,y_i,y_j)\mid$\

\hspace{4.8cm} $ \geq s(a_1+1)- {s \choose 2}(c_2-1)$.

\end{pf}

   \begin{pro} \label{2}
      \rm{Let} $\Gamma$ be a distance-regular graph with valency $k$ and diameter $D \geq 2$. Let $s$ be maximal such that for all $x$ and all $y$, $z$ $\in \Gamma(x)$ with $y\nsim z$, there exists a co-clique of size at least $s$ in $\Delta(x)$ containing $y$ and $z$. Then
      \begin{description}
        \item [$(i)$] $s \geq \frac{k}{a_1+1}$
        \item [$(ii)$] $c_2-1 \geq$ max$\{\frac{s^{'}(a_1+1)-k}{{s^{'} \choose 2}}\mid 2 \leq s^{'} \leq s \}$ and equality implies $\Gamma$ is a Terwilliger graph.
      \end{description}
   \end{pro}
    \begin{pf}
     Let $s$ be maximal satisfying the condition in the Proposition~\ref{2}. Then $k \leq s(a_1+1)$ as $\Delta(x)$ has valency $a_1$ and $k$ vertices. This shows $(i)$. By Lemma~\ref{1}, the inequality in $(ii)$ holds. Next, we assume that the equality holds in $(ii)$.  Let $2\leq s^{''}\leq s$ be an integer satisfying $\frac{s^{''}(a_1+1)-k}{{s^{''} \choose 2}}$ = max$\{\frac{s^{'}(a_1+1)-k}{{s^{'} \choose 2}}\mid 2 \leq s^{'} \leq s \}$, then by Lemma~\ref{1} there exists a co-clique $\bar{C}^{''}$ on $\{y_1, y_2, \cdots, y_{s^{''}}\}$ such that for any two vertices $y_i, y_j$ at distance two, $\mid\Gamma(x,y_i,y_j)\mid$ = $c_2-1$ holds. That is, if we take three vertices $z_1, z_2$ and $z_3$ such that $d_{\Gamma}(z_2,z_3)=2$, $z_1\sim z_2$ and $z_1\sim z_3$, then since $z_1\in \Gamma(z_2,z_3)$ and $\mid\Gamma(z_1,z_2,z_3)\mid$ = $c_2-1$, the valency of $z_1$ in $\Gamma(z_2,z_3)$ is $c_2-1$. Hence the subgraph induced by $\Gamma(z,w)$ is a clique of size $c_2$ for any two vertices $z$ and $w$ at distance two in $\Gamma$. So, $\Gamma$ is a Terwilliger graph.
    \end{pf}

    For all known examples of Terwilliger graphs, we have equality in case $(ii)$ above.

\begin{thm} \label{3}
 \rm{Let} $\Gamma$ be a distance-regular graph with valency $k$ and diameter $D \geq 2$ and define $\alpha = \lceil \frac{k}{a_1+1}\rceil$. Then  $c_2-1 \geq \frac{\alpha(a_1+1)-k}{{\alpha \choose 2}}$ and equality implies that $\Gamma$ is a Terwilliger graph.
\end{thm}
\begin{pf}
This is an immediate consequence of Proposition~\ref{2}.
\end{pf}

 Theorem~\ref{3} gives the new existence condition for distance-regular graphs and the following two intersection arrays in \cite[p.425-431]{bcn} are ruled out.

\begin{cor}
\rm{There} are no distance-regular graphs with one of the following intersection arrays:\\
\begin{center}
$(i)\hspace{0.2cm} \{44,30,5;1,3,40\}$; ~~~~~~~~~~~~~   $(ii)\hspace{0.2cm} \{65,44,11;1,4,55\}$.
\end{center}
\end{cor}
\begin{pf}
\begin{enumerate}
  \item Since a distance-regular graph $\Gamma$ with intersection array $(i)$ satisfies $c_2-1 = 2 = \frac{\alpha(a_1+1)-k}{{\alpha \choose 2}}$, where $\alpha={\lceil \frac{k}{a_1+1}\rceil}=4$, $\Gamma$ is a Tewilliger graph by Theorem~\ref{3}. But this is impossible by \cite[Corollary 1.16.6]{bcn}.
  \item As $\lceil \frac{65}{20+1}\rceil$ = 4, there is no distance-regular graph with intersection array $(ii)$, by Theorem~\ref{3}.
\end{enumerate}

\end{pf}

\hspace{-0.5cm}\textbf{Remark}: A. Juri$\check{\rm{s}}$i$\acute{\rm{c}}$ and J. Koolen \cite{kj} proved that distance-regular graphs with intersection arrays \{81,56,24,1;1,3,56,81\}, \{117,80,30,1;1,6,80,117\}, \{117,80,32,1;\\1,4,80,117\}, and \{189,128,45,1;1,9,128,189\} do not exist. It also follows from Theorem~\ref{3}.

\section{Shilla distance-regular graphs}

~~~~~In this section we first give a lower bound on the second largest eigenvalue of a distance-regular graph with diameter 3. Then we will define Shilla distance-regular graphs and give some results on them.

\begin{lemma}\label{10}
 \rm{Let} $\Gamma$ be a distance-regular graph with valency $k$ and diameter 3. Then the second largest eigenvalue $\theta_1$ of $\Gamma$ satisfies
 \begin{center}
 $\theta_{1}\geq$ min$\{\frac{a_{1}+\sqrt{a_{1}^{2}+4k}}{2}, a_3\}$.\\
\end{center}
\end{lemma}
 \begin{pf}
 Let $x$ be a vertex of $\Gamma$. As the induced subgraph on $\{x\} \cup \Gamma(x)$ (respectively $\Gamma_3(x)$) has largest eigenvalue $\frac{a_{1}+\sqrt{a_{1}^{2}+4k}}{2}$ (respectively $a_3$), it follows that the induced subgraph on $\{x\} \cup \Gamma(x) \cup \Gamma_3(x)$ has the second largest eigenvalue at least min$\{\frac{a_{1}+\sqrt{a_{1}^{2}+4k}}{2}, a_3\}$. Now the lemma follows by Theorem~\ref{0}.
 \end{pf}\\
\vspace{-1cm}
\begin{theorem}\label{11}
\rm Let $\theta$ be an eigenvalue of a distance-regular graph $\Gamma$ with valency $k$ and diameter 3. Then the following are equivalent:
\begin{description}
  \item[$(i)$] $u_2(\theta)=0$;
  \item[$(ii)$] $\theta=a_3$;
  \item[$(iii)$] $\theta=\frac{a_{1}+\sqrt{a_{1}^{2}+4k}}{2}$.
\end{description}
Moreover, if $(i)-(iii)$ hold then $\theta$ is the second largest eigenvalue of $\Gamma$.
\end{theorem}
\begin{pf}
$(i)\Leftrightarrow(ii)$: Let $u_0=1$, $u_1$, $u_2$, $u_3$ be the standard sequence for $\theta$. As $c_3u_2+a_3u_3=\theta u_3$, it follows $u_2=0$ if and only if $\theta=a_3$.\\
$(i)\Leftrightarrow(iii)$: As $1+a_1u_1+b_1u_2=\theta u_1$
and $u_1=\frac{\theta}{k}$, it follows that $\theta=\frac{a_{1}+\sqrt{a_{1}^{2}+4k}}{2}$ if and only if $u_2=0$.\\
Moreover, if $(i)-(iii)$ hold then, $u_0=1>0$, $u_1=\frac{a_3}{k}>0$, $u_2=0$, and $u_3=-\frac{c_2a_3}{kb_2}<0$ and hence $\theta$ is the second largest eigenvalue of $\Gamma$ by \cite[Corollary 4.1.2]{bcn}.

\end{pf}

A distance-regular graph with diameter 3 and valency $k$ is called {\em Shilla} if its second largest eigenvalue $\theta_1$ satisfies $\theta_1=a_3$. It follows by Theorem~\ref{11} that $\theta_1 =a_3=\frac{a_{1}+\sqrt{a_{1}^{2}+4k}}{2}$ and hence $k=(a_3-a_1)a_3$. For a Shilla distance-regular graph $\Gamma$, put $b=b(\Gamma):=a_3-a_1$. Then clearly $b\geq2$ and $k=ba_3$.\\

A Shilla distance-regular graph $\Gamma$ with $b(\Gamma)=b$ has distinct four eigenvalues $\theta_0=k=ba_3>\theta_1=a_3>\theta_2>\theta_3$, where $\theta_2$ and $\theta_3$ are  two roots of the equation $x^2-(a_1+a_2-k)x+(b-1)b_2-a_2=0$. Let $m_i$ be the multiplicity of $\theta_i$.\\
 If both $\theta_2$ and $\theta_3$ are integers then  $(a_1+a_2-k)^2-4((b-1)b_2-a_2)$ is a perfect square.\\ If both $\theta_2$ and $\theta_3$ are  non-integers, then $m_2=m_3$ holds. This implies, by Equation~\ref{e0}, that the equation\\
 \begin{equation}
 \begin{split}
 &~(b_2+c_2)(b_2+c_2-a_3)(b_2+c_2+(b-1)a_3)-bb_2^2+(2b-3)c_2^2\\
 & +b(b-1)c_2+(b-1)^2a_3c_2-b(b-1)a_3b_2+(b-3)b_2c_2=0
 \end{split}
 \label{e1}
 \end{equation}
 holds. In Theorem~\ref{13} below, we will discuss the situation $m_2=m_3$ in more detail.\\

 Now, we will show that there are finitely many Shilla distance-regular graphs $\Gamma$ with fixed $b(\Gamma)$. To do this we first show Lemma~\ref{4}.

\begin{lemma}\label{4}
\rm{Let} $\Gamma$ be a Shilla distance-regular graph with $b(\Gamma)=b$. Then,\\ $$c_2 \geq \frac{2a_3-b^2+b+2}{b(b+1)}.$$
\end{lemma}
\begin{pf}
Let $x$ be a vertex of $\Gamma$. Then there exists a co-clique of size $b+1$ in $\Delta(x)$ as $k=ba_3=b(a_1+b)>b(a_1+1)$ and by Lemma~\ref{1}, the proof is complete.
\end{pf}

\begin{thm}\label{thm}
\rm{For} given $\beta \geq 2$, there are finitely many Shilla distance-regular graphs $\Gamma$ with $b(\Gamma)=\beta$.
\end{thm}
\begin{pf}
For given $\beta \geq 2$, let $\Gamma$ be a Shilla distance-regular graph with valency $k$, $b(\Gamma)=\beta$ and $n$ vertices. Then clearly $k=\beta a_3=\beta(a_1+\beta)=\beta(a_1+1)+{\beta}^2-\beta$. We will show that $k$ is bounded above by $\beta$. We first show the following.\\
\begin{center}
\hspace{-4.9cm}\textbf{Claim}:\hspace{1.5cm} $k<\beta^3-\beta$~~~~ or~~~~ $n<k(2\beta^3-\beta+1)$.
\end{center}

 \hspace{-0.6cm}\textbf{Proof of claim}: If $a_1+1<\beta^2-\beta $, then $k=\beta a_3=\beta(a_1+\beta)<{\beta}^3-\beta$. So, let us assume $ a_1+1\geq{\beta}^2-\beta $ then clearly $k\geq{\beta}^3-\beta $. Lemma~\ref{4} implies $c_2 \geq \frac{a_3+(a_3+1-{\beta}^2)+\beta+1}{\beta(\beta+1)}>\frac{a_3+1}{\beta(\beta+1)}$, where the second inequality follows from $a_3+1\geq {\beta}^2$. As $c_3=(\beta-1)a_3$ and $b_1=(\beta-1)(a_3+1)$, it follows that
\begin{center}
\hspace{-4.1cm}$n=1+k+k\frac{b_1}{c_2}+k\frac{b_1b_2}{c_2c_3}=1+k+k\frac{(\beta-1)(a_3+1)}{c_2}+k\frac{b_2}{c_2}+\frac{\beta b_2}{c_2}$\\
$\leq 1+k+2k\beta(\beta-1)(\beta+1)+{\beta}^2({\beta}^2-1)\leq 1+k+2k\beta(\beta-1)(\beta+1)+k\beta$\\
\hspace{-10.7cm}$< k(2{\beta}^3-\beta+2).$
\end{center}
So, the claim is proved.\\

Now by letting $m_1$ to be the multiplicity of $\theta_1=a_3$, it follows from Equation~\ref{e0} that $m_1<\frac{n}{u_1(a_3)^2k}$. By \cite[Theorem 5.3.2]{bcn}, $\sqrt{k} < m_1$ holds. As $m_1<\frac{n}{u_1(a_3)^2k}\\ <\frac{k(2{\beta}^3-\beta+2)}{u_1(a_3)^2k}=2{\beta}^5-{\beta}^3+2{\beta}^2$, it follows $k < 4{\beta}^{10}$. This shows the theorem.
\end{pf}

In the next result, we give some divisibility conditions for Shilla distance-regular graphs.

 \begin{lemma}\label{12}
     \rm{Let} $\Gamma$ be a Shilla distance-regular graph with $b(\Gamma)=b$. Then the following holds:
     \begin{description}
       \item $(i)$~~  $c_2$ divides $(b-1)a_3b_2$;
       \item $(ii)$~ $c_2$ divides $(b-1)ba_3(a_3+1)$;
       \item $(iii)$~$c_2$ divides  $b(a_3+1)b_2$;
       \item $(iv)$~ $c_2$ divides $(b+a_3)b_2$ and $(b+a_3)b_2 \geq (1+a_3)c_2$, where equality is attained if and only if $p^3_{33}=0$;
       \item $(v)$~~ $c_2$ divides $(b-1)bb_2$.
     \end{description}
     \end{lemma}

     \begin{pf}
     $(i)$ follows from the fact that $p^3_{32}$ is a non-negative integral.

     $(ii)$ and $(iii)$ hold as $k_2=\frac{kb_1}{c_2}$ and $k_3=\frac{kb_1b_2}{c_2c_3}$ are integral respectively.
    Since $\frac{kb_1b_2}{c_2c_3}$\\$=k_3=p^3_{30}+p^3_{31}+p^3_{32}+p^3_{33}=1+a_3+\frac{c_3(b_2-1)+a_3(a_3-1-a_1)}{c_2}+p^3_{33} $ is integral, it follows that $\frac{(b+a_3)b_2}{c_2} \geq 1+a_3$ and equality is attained if and only if $p^3_{33}=0$. Since $\frac{(b+a_3)b_2}{c_2}$ and $\frac{(b-1)(b+a_3)b_2}{c_2}$ are integral, $(v)$ holds by $(i)$.
     \end{pf}\\

\vspace{-0.5cm}
We will give some necessary conditions for Shilla distance-regular graphs with $m_2=m_3$.

\begin{thm}\label{13}
\rm{Let} $\Gamma$ be a Shilla distance-regular graph with $b(\Gamma)=b\geq2$ and\\ $m_2=m_3$. Then the following holds:\\
\vspace{-0.5cm}
\begin{description}
  \item $(i)$ $a_3-b<b_2+c_2<a_3+b;$
  \vspace{-0.2cm}
  \item $(ii)$ $ c_2< b_2+b;$
  \vspace{-0.2cm}
  \item $(iii)$ If $b_2+c_2=a_3$ or $b_2=c_2$, then $a_3=\frac{b(b-1)}{2}$.
\end{description}

\end{thm}
\begin{pf}
Since $m_2=m_3$, Equation(\ref{e1}) holds. \\
$(i)$: If $b_2+c_2\leq a_3-b$, then the LHS of Equation(\ref{e1}) is negative. Hence\ $a_3-b<b_2+c_2$. In similar fashion, we can show $b_2+c_2<a_3+b$. Thus, $a_3-b<b_2+c_2$ \ $<a_3+b$.\\
$(ii)$:If $c_2\geq b_2+b$, then Lemma~\ref{12} $(iv)$ implies $(b+a_3)b_2\geq (1+a_3)c_2\geq(1+a_3)(b_2+b)$. This means, $b_2\geq \frac{b}{b-1}(a_3+1)=a_3+1+\frac{a_3+1}{b-1}>a_3+2,$ where the last inequality holds by $a_3\geq b$. So, $b_2+c_2>a_3+b$ and this is a contradiction to $(i)$. Thus $c_2<b_2+b$.\\
$(iii)$ If $b_2+c_2=a_3$, then Equation(\ref{e1}) becomes $b^2(b_2^2-c_2^2)+2c_2(b_2+c_2)=b(b-1)c_2$ and it follows $b_2\leq c_2$. If $b_2=c_2$, then $c_2=\frac{b(b-1)}{4}$ and hence $a_3=\frac{b(b-1)}{2}$. If $b_2<c_2$, then $b(b-1)c_2=b^2(b_2^2-c_2^2)+2c_2(b_2+c_2)\leq 4c_2^2-2c_2-2b^2c_2+b^2$ and hence $c_2\geq\frac{2b^2-b+2}{4}$. Now it follows from Lemma~\ref{12} that $\frac{(b+a_3)b_2}{c_2}=\frac{(b+b_2+c_2)b_2}{c_2}$ is integral and hence  $\frac{(b_2+b)b_2}{c_2}$ is integral. Since $b_2=c_2-\alpha$ for some $1\leq \alpha< b$, we find $c_2$ divides $\alpha(b-\alpha)$. Hence $c_2\leq\alpha(b-\alpha)\leq\frac{b^2}{4}$, but this contradicts $c_2\geq\frac{2b^2-b+2}{4}.$
If $b_2=c_2$, then Equation(\ref{e1}) becomes $2c_2(2c_2-a_3)(2c_2+(b-1)a_3)+b(b-1)c_2-(b-1)a_3c_2+2(b-3)c_2^2$ and this is always positive (respectively negative) if $2c_2>a_3$ (respectively $2c_2<a_3$). Thus $2c_2=a_3$ and hence $a_3=\frac{b(b-1)}{2}$.
\end{pf}

Note that in case $(iii)$ above, we have the intersection arrays \\$$\{\frac{b^2(b-1)}{2},\frac{(b-1)(b^2-b+2)}{2},\frac{b(b-1)}{4};1,\frac{b(b-1)}{4},\frac{b(b-1)^2}{2}\}$$  and they are only feasible for $b\equiv0,1$ $(\mod4)$. Besides this family of intersection arrays, using the computer, the only  other feasible intersection arrays for Shilla distance-regular graphs with $m_2=m_3$ and $a_3\leq 100$ are:

\begin{center}
$\begin{array}{ccc}
  \hspace{-1.5cm}(\divideontimes)& \hspace{-0.5cm}(i)~~ \{120,117,20;1,1,108\}; & \hspace{1cm}(ii)~~ \{676,675,31;1,9,650\}; \\
 &\hspace{-0.05cm}(iii)~~ \{486,440,50;1,10,432\}; & \hspace{2.1cm}(iv)~~ \{4264,4233,102;1,17,4182\}. \\
\end{array}$
\end{center}

In the next theorem, we classify  Shilla distance-regular graphs $\Gamma$ with $b(\Gamma)=2$.

\vspace{0.3cm}

\begin{thm}\label{20}
\rm Let $\Gamma$ be a Shilla distance-regular graph with $b(\Gamma)=2$. Then $\Gamma$ is one of the following graphs:\\
\vspace{-0.5cm}
\begin{description}
    \item [$(i)$] the Odd graph with valency 4;
    \vspace{-0.2cm}
    \item [$(ii)$] a generalized hexagon of order (2,2);
    \vspace{-0.2cm}
    \item [$(iii)$] the Hamming graph $H(3,3)$;
    \vspace{-0.2cm}
    \item [$(iv)$] the Doro graph with intersection array $\{10,6,4;1,2,5\}$;
    \vspace{-0.2cm}
    \item [$(v)$] the Johnson graph $J(9,3)$.
\end{description}
\end{thm}

    \begin{pf}
    If $b=2$, then $b_1=a_3+1$ and hence $\theta_1=a_3=b_1-1$. By \cite[Theorem 4.4.11]{bcn} we only need to consider the cases $c_2=1$ or $\Gamma$ is one of a Doob, a Hamming, a locally Petersen, a Johnson, a halved cube, or the Gosset graph. By \cite[Theorem 1.16.5]{bcn} the only locally Petersen graph which is a Shilla distance-regular graph is the Doro graph. If $\Gamma$ is a Doob or a Hamming graph, then $c_3=3=a_3$ and it follows that $\Gamma=H(3,3)$. If $\Gamma$ is a Johnson graph, then $c_3=9=a_3$ and it follows that $\Gamma=J(9,3)$. Neither the Gosset graph nor the halved 6-cube are possible as they have $a_3=0$. Also, the halved 7-cube is not a Shilla distance-regular graph. To complete the proof of this theorem, we only need to consider the case $c_2=1$. If $c_2=1$, then $\Gamma$ is a locally disjoint union of $(a_1+1)$-cliques. This implies that $a_1+1$ divides $k$, and hence $a_3 \in \{2,3\}$. If $a_3=2$, then $k=4, c_1=1, a_1=0, b_1=3, c_3=2$ and $b_2\in \{1,2,3 \}.$ Only for $b_2=3$, the multiplicity $m_1$ is an integer, and $\Gamma$ is the Odd graph with valency 4. If $a_3=3$, then $k=6, c_1=1, a_1=1, b_1=4 $, and $b_2 \in \{1,2,3,4,5 \}$. Only for $b_2=4$, the multiplicity $m_1$ is an integer, and $\Gamma$ is a generalized hexagon of order (2,2).
   \end{pf}\\

\vspace{-0.5cm}

The following lemma gives a sufficient condition for a distance-regular graph to be $Q$-polynomial.
 \begin{lemma}\label{16}
 \rm{Let} $\Gamma$ be a distance-regular graph with diameter $D=3$, $n$ vertices and eigenvalues $k>\theta_1>\theta_2>\theta_3$. If $n\geq \frac{(m_1+2)(m_1+1)}{2}$, then $q^2_{11}=0$ or $q^3_{11}=0$ and hence $\Gamma$ is $Q$-polynomial with respect to $\theta_1$, where $q^2_{11}$ and $q^3_{11}$ are Krein parameters.
\end{lemma}
    \begin{pf}
    As $n\geq \frac{(m_1+2)(m_1+1)}{2}$ and $n=1+m_1+m_2+m_3$, it follows that $m_2+m_3$ \\$\geq{m_1+1\choose2}.$  As $\displaystyle \sum_{ q^i_{11}\neq 0} m_i\leq $ $m_1+1 \choose 2$~\cite[Proposition 4.1.5]{bcn} and $q^0_{11}>0$, it follows that $q^2_{11}=0$ or $q^3_{11}=0$. This implies that $\Gamma$ is $Q$-polynomial with respect to $\theta_1$.
    \end{pf}\\
\vspace{-0.5cm}
\begin{lemma}\label{17}
\rm Let $\Gamma$ be a Shilla distance-regular graph with $n$ vertices, $b(\Gamma)=\beta$ and valency $k$. If $\Gamma$ is not $Q$-polynomial with respect to $\theta_1$ then $k<\beta^5(\beta+1)^2 $.
\end{lemma}

\begin{pf}
Let $k\geq \beta^3-\beta$. If $\Gamma$ is not $Q$-polynomial with respect to $\theta_1$ then $n<\frac{(m_1+1)(m_1+2)}{2}$ by Lemma~\ref{16}. By Equation~\ref{e0} we find $m_1+2<\frac{n}{a_3/\beta}+2<\frac{\beta+1}{\beta}\frac{n}{a_3/\beta}$, where the last inequality holds by $\frac{n}{a_3}>\frac{k}{a_3}=\beta\geq 2$. Combining the above two inequalities we find $\sqrt{2n}<\frac{\beta+1}{\beta}\frac{n}{a_3/\beta}$ and hence $2(\frac{a_3}{\beta+1})^2<n$. As $k\geq \beta^3-\beta$, by Claim in Theorem~\ref{thm}, we find $n<k(2\beta^3-\beta+1)$. So $2(\frac{a_3}{\beta+1})^2<k(2\beta^3-\beta+1)$. Since $k=\beta a_3$, we find $k<\beta^5(\beta+1)^2 $.
\end{pf}

    \begin{pro} \label{pro}
    \rm{Let} $\Gamma$ be a Shilla distance-regular graph with $b(\Gamma)=b$. Then the following holds:
  \begin{description}
    \item[$(i)$] $q^2_{11}>0$;
    \item[$(ii)$] $q^3_{11}\geq 0$ if and only if $\theta_3 \geq -\frac{b(bb_2+c_2)}{b_2+c_2}$.\
  \end{description}

  So, in particular $\theta_3 \geq -\frac{b(bb_2+c_2)}{b_2+c_2}$ holds.
  \end{pro}
     \begin{pf}
     Note that $q^i_{jh}\geq 0=\frac{m_jm_h}{|V\Gamma|}\displaystyle \sum_{l=0}^D k_lu_l(\theta_i)u_l(\theta_j)u_l(\theta_h)\geq 0$~\cite[Proposition 4.1.5]{bcn}. Hence $q^i_{11}\geq 0$ if and only if $\displaystyle \sum_{l=0}^3 k_lu_l(\theta_1)^2u_l(\theta_i)\geq 0$ if and only if \
     $$c_2\theta_i^3-c_2(a_1+a_2)\theta_i^2+(c_2a_1a_2-b_1c_2^2-c_2kb_2^2c_3)\theta_i+c_2a_2k+b^2b_2^2c_3\geq0$$

     Since $\theta_2$ and $\theta_3$ are two roots of polynomial $\theta^2-(a_1+a_2-k)\theta+(b-1)b_2-a_2$, we obtain that $q^i_{11}\geq 0$ if and only if $(b^2b_2+bc_2+b_2\theta_i+c_2\theta_i)\geq0$ for $i=2,3.$ As, $\theta_2>\theta_3$ we see immediately that $q^2_{11}>0$ and also that $q^3_{11}\geq 0$ if and only if $\theta_3 \geq -\frac{b(bb_2+c_2)}{b_2+c_2}$. This shows the proposition.
      \end{pf}\\
\vspace{-0.7cm}
 \begin{cor}\label{14}
   \rm{Let} $\Gamma$ be a Shilla distance-regular graph with $b(\Gamma)=b$. Then \\
   $$-b^2 < \theta_3 < -b$$
     \end{cor}
     \begin{pf}
     Let $x$ be a vertex in $\Gamma$. Then the induced subgraph on $\{x\} \cup \Gamma(x)$ has two eigenvalues, $-b$ and $a_3$. Then, by Theorem~\ref{0}, $ \theta_3\leq -b $ holds. But if $\theta_3=-b$, then $u_2(\theta_3)=0$ and this is not possible by Theorem~\ref{11}. The lower bound follows immediately from Proposition~\ref{pro}.
     \end{pf}\\

We will improve the lower bound for $\theta_3$ in Theorem~\ref{15} below.
   \begin{cor}\label{8}
   \rm{Let} $\Gamma$ be a Shilla distance-regular graph with $b(\Gamma)=b$. Then $\Gamma$  \rm{is}   $Q$-polynomial with respect to $\theta_1$  if and only if  $\theta_3=-\frac{b(bb_2+c_2)}{b_2+c_2}$. If $\Gamma$ is  $Q$-polynomial with respect to $\theta_1$, then all eigenvalues of $\Gamma$ are integral, $b_2+c_2$ divides $b(b-1)b_2$ and $-b^2+1\leq\theta_3\leq-\frac{b^2(b+3)}{3b+1}$.\\
   \end{cor}
   \begin{pf}
   Note that $\Gamma$ is  $Q$-polynomial with respect to $\theta_1$ if and only if $q^2_{11}=0$ or $q^3_{11}=0$.  Hence the first part of this corollary follows from Proposition~\ref{pro}. Assume $\Gamma$ is $Q$-polynomial with respect to $\theta_1$. Then $\theta_3=-b-\frac{b(b-1)b_2}{b_2+c_2}$ and hence $\theta_3$ is integral. So, all the eigenvalues are integral. As $b\leq a_3$, it follows that $c_2\leq\frac{(b+a_3)b_2}{1+a_3}\leq\frac{2bb_2}{1+b}$ by Lemma~\ref{12}. So, $\theta_3\leq -b-\frac{b(b-1)(b+1)}{3b+1}= -\frac{b^2(b+3)}{3b+1}$. Thus $-b^2+1\leq\theta_3\leq-\frac{b^2(b+3)}{3b+1}$ holds by Corollary~\ref{14}.
  \end{pf}

In the next two results, we classify the Shilla distance-regular graphs $\Gamma$ with $b(\Gamma)=3$.

\begin{pro}\label{pro2}
 Let $\Gamma$ be a Shilla distance-regular graph with $b(\Gamma)=3$ and let $\Gamma$ be  $Q$-polynomial with respect to $\theta_1$. Then $\Gamma$ has one of the following intersection arrays.\\
 \begin{center}
 $\begin{array}{ccc}
    (i)~~ \{ 42, 30, 12; 1, 6, 28 \},\hspace{0.5cm} & (ii)~~ \{ 105, 72, 24; 1, 12, 70 \}.
 \end{array}$
 \end{center}

\end{pro}
 \begin{pf}
 By Corollary~\ref{8}, if $\Gamma$ is a $Q$-polynomial with respect to $\theta_1$, then $\theta_3\in\{ -6, -7, -8\}$.\

   If $\theta_3=-6$, then $c_2=b_2$. Since $\theta_3$ is a root of the equation $ x^2-(a_1+a_2-k)x$\ $+(b-1)b_2-a_2=0$, $b_2=c_2$ and $b(\Gamma)=b=3$, it follows $a_3=\frac{8}{3}b_2-6$. But this means $m_1=33-\frac{135}{2b_2}$ what is impossible, as $m_1$ must be an integer. \

 If $\theta_3=-7$,  Similarly, we obtain $b_2=2c_2$ and $a_3=\frac{7}{2}c_2-7$. Then $m_1=60-\frac{108}{c_2}$. Since $a_3$ and $m_1$ have to be integers, it follows that $2$ divides $c_2$ and $c_2$ divides $108$, that is, $c_2 \in \{2, 4, 6, 12, 18, 36, 54, 108\}$. The case $c_2=2$ implies $a_3=0$, which is impossible. For $c_2\in\{18, 36, 54, 108\}$, we find that $m_3$ is non-integral.   The case $c_2=4$ gives us the intersection array $\{ 21, 16, 8; 1, 4, 14 \}$ and it was shown by K. Coolsaet\cite{cool} that a distance-regular graph with this intersection array does not exist. The cases $c_2=6$ and $c_2=12$ give the intersection arrays $(i)$ and $(ii)$ respectively.\

If $\theta_3=-8$,  Similarly, we obtain $b_2=5c_2$ and $a_3=\frac{32}{5}c_2-8$. But this means $m_1=141-\frac{315}{2c_2}$ what is impossible.\\
 \end{pf}\\

\vspace{-1cm}
 \begin{thm}\label{21}
      \rm{Let} $\Gamma$ be a Shilla distance-regular graph with $b(\Gamma)=3$. Then $\Gamma$ has  one of the following intersection arrays.\\

      \begin{center}

        \hspace{-0.4cm}$(i)\hspace{0.2cm} \{12, 10, 5; 1, 1, 8\},$ ~~~~~~ $(ii)\hspace{0.2cm} \{12, 10, 2; 1, 2, 8 \},$~~~~~~ $(iii)\hspace{0.2cm} \{12, 10, 3; 1, 3, 8\},$ \\

        \hspace{-0.05cm}$(iv)\hspace{0.2cm} \{15, 12, 6; 1, 2, 10\},$~~~~~\hspace{-0.15cm} $(v)\hspace{0.2cm} \{24, 18, 9; 1, 1, 16\}, $ ~~~\hspace{0.3cm}$(vi)\hspace{0.2cm} \{27, 20, 10; 1, 2, 18\}$,  \\

        \hspace{-0.25cm}$(vii)\hspace{0.2cm} \{30, 22, 9; 1, 3, 20\},$~~~ $(viii)\hspace{0.2cm} \{42, 30, 12; 1, 6, 28\},$ \hspace{0.05cm}$(ix) \{60, 42, 18; 1, 6, 40\}$,  \\

        \hspace{0.5cm}$(x)\hspace{0.2cm} \{69, 48, 24; 1, 4, 46\},$~~~ $(xi)\hspace{0.2cm} \{93, 64, 24; 1, 6, 62\},$~~~ $(xii)\hspace{0.2cm}  \{105, 72, 24; 1, 12, 70\}.$ \\

      \end{center}
    \end{thm}

Note that all the above intersection arrays have $\theta_3\geq-7$\\

\begin{pf}
If $\Gamma$ is $Q$-polynomial with respect to $\theta_1$ then it follows from Proposition~\ref{pro2} that $\Gamma$ has intersection array $(viii)$ or $(xii)$. If $\Gamma$ is not $Q$-polynomial with respect to $\theta_1$ then by Lemma~\ref{17}, $a_3<3^4\times4^2=1296$. We checked by computer that the above arrays are the only possible intersection arrays for Shilla distance-regular graphs with $a_3<1296$.
\end{pf}

\textbf{Remark}: The unitary nonisotropics graph with $q=4$ as defined in \cite[Section 12.4]{bcn} has intersection array $(i)$. It is not know whether it is unique or not. There exists a unique distance-regular graph with intersection array $(iii)$ namely the Doro graph as defined in \cite[Section 12.1]{bcn}. For the other intersection arrays, it is not known whether a distance-regular graph with those intersection arrays does exist, or not.\\

Now, we improve the lower bound of the smallest eigenvalue $\theta_3$ for a Shilla distance-regular graph.

\begin{thm}\label{15}
\rm{For} a Shilla distance-regular graph with $b(\Gamma)=b$ and smallest eigenvalue $\theta_3$, we have $\theta_3<-b^2+2$ if and only if $b=2$.
\end{thm}
\begin{pf}
($\Leftarrow$) For $b=2$, we are done by Theorem~\ref{20}.

($\Rightarrow$) Let $\theta_3<-b^2+2$. Then by Theorems~\ref{20} and \ref{21} we have $b=2$ or $b\geq 4$. So, let us assume $b\geq 4$. By Corollary~\ref{14}, we have $-b^2<\theta_3<-b^2+2$. Then either $\theta_3=-b^2+1$ or $m_2=m_3$ and $\theta_3$ is non-integer. If $\theta_3=-b^2+1$, then by Proposition~\ref{pro}, $b_2\geq(b^2-b-1)c_2$. Since $\theta_3=-b^2+1$ is a root of the equation $ x^2-(a_1+a_2-k)x$\ $+(b-1)b_2-a_2=0$ and $b_2\geq(b^2-b-1)c_2$, we have\\  $a_3=b_2+\frac{b^2-2}{b^2-b-1}c_2-(b^2-1)\geq\frac{(b-1)^3(b+1)}{(b^2-b-1)}c_2-(b^2-1)$ \\and hence $c_2\in\{1,2,3\}$ by Lemma~\ref{12}. Since $a_3=b_2+\frac{b^2-2}{b^2-b-1}c_2-(b^2-1)\in\textbf{Z}$, $b^2-b-1$ divides $(b-1)c_2$ and hence $b^2-b-1\leq 3(b-1)$. This contradicts $b\geq 4$. Let us now consider the case $m_2=m_3$ and $\theta_3$ is non-integer. As $c_2+b_2\leq a_3+b-1$, by Theorem~\ref{13}, the LHS of Equation~(\ref{e1}) is at most
\vspace{0cm}
\begin{center}
$\begin{array}{cc}
(\divideontimes\divideontimes)&\hspace{0.8cm}(b_2+c_2)((3b-4)c_2-b_2)+(b-1)((2b-2)c_2-b_2)a_3+b(b-1)c_2.\\
\end{array}$
\end{center}
Since $b_2\geq(3b-4)c_2$ implies that $(\divideontimes\divideontimes)$ is negative, we have $b_2<(3b-4)c_2$. By Proposition~\ref{pro}, we have $\theta_3\geq -b^2+\frac{b(b-1)c_2}{b_2+c_2}\geq -b^2+\frac{b}{3}$, and hence $4\leq b\leq 5$. Let us consider first $b=5$ . As $-23>\theta_3\geq -25+\frac{20c_2}{b_2+c_2}$, it follows that $9c_2<b_2< 11c_2$. As $(\divideontimes\divideontimes)>0$, it follows that $a_3\leq 13$. As the intersection array $\{50,44,5;1,5,40\}$ has $\theta_3=-7.623$, this case follows now from $(\divideontimes)$.  Secondly, we assume that $b=4$. As $-14>\theta_3\geq-16+\frac{12c_2}{b_2+c_2}$, it follows that $5c_2< b_2< 8c_2$. Since $b_2+c_2\leq a_3$ implies that the LHS of Equation~(\ref{e1}) is negative, we have $a_3<b_2+c_2\leq a_3+3$. If $b_2+c_2=a_3+1$ then the LHS of Equation~(\ref{e1})=$4a_3^2+5a_3+1+9a_3c_2-12a_3b_2-4b_2^2+5c_2^2+12c_2+b_2c_2$. Since $5c_2\leq b_2$, $b_2^2\geq 5c_2^2+12c_2+b_2c_2$ and $6b_2\geq 5a_3 $. Hence the LHS of Equation~(\ref{e1}) is negative. Similarly, for $b_2+c_2=a_3+2$, the LHS of Equation~(\ref{e1}) is negative if $a_3\geq 9$. So, $a_3\leq 8$ and then we are done by $(\divideontimes)$. If $b_2+c_2=a_3+3$ then the LHS of Equation~(\ref{e1})=$2(13c_2-2b_2)(c_2+b_2)+3(3b_2-5c_2-9)$. Since $b_2\leq\frac{19}{3}c_2$ implies that the LHS of Equation~(\ref{e1})is positive whenever $b_2+c_2\geq 45$ (i.e.$a_3\geq42$), either ($b_2\leq\frac{19}{3}c_2$ and $a_3\leq 41$) or $b_2>\frac{19}{3}c_2$. By  $(\divideontimes)$ we obtain $b_2>\frac{19}{3}c_2$. Since $b_2\geq\frac{27}{4}c_2$ implies that the LHS of Equation~(\ref{e1}) is negative whenever $c_2\geq 9$, either ($b_2\geq \frac{27}{4}c_2$ and $c_2\leq 8$) or $b_2<\frac{27}{4}c_2$. As $c_2\leq 8$ implies that $a_3\leq 85$ by Lemma~\ref{4}, from $(\divideontimes)$  we obtain $\frac{19}{3}c_2<b_2<\frac{27}{4}c_2$. Then by Lemma~\ref{12} $(v)$, $\frac{12b_2}{c_2}\in\{77,78,79,80\}$. In the first two possibilities the LHS of Equation~(\ref{e1}) is negative. In the last two possibilities the number $c_2$ is non-integral. This shows the theorem.
\end{pf}

\begin{flushleft}
\Large\textbf{Acknowledgments}
\end{flushleft}

Section 3 was inspired by a conversation with Sejeong Bang. We also would like to thank her for the careful reading she did. Comments by Mitsugu Hirasaka are appreciated much. This work is partially supported by KRF-2008-314-C00007 and this support is greatly appreciated.

\clearpage

\end{document}